\documentclass[12pt]{amsart}
\usepackage{euscript}
\usepackage{amssymb,amscd}

\newtheorem{theorem}{Theorem}[section]

\newtheorem{proposition}[theorem]{Proposition}

\begin{document}
\title[Thermodynamical Formalism for Inducing Schemes]
{Thermodynamical Formalism Associated with Inducing Schemes For
  One-dimensional Maps}

\author{Yakov Pesin}
\address{Department of Mathematics, Pennsylvania State
University, University Park, State College, PA 16802}
\email{pesin@math.psu.edu}
\author{Samuel Senti}
\address{Instituto de Matematica Pura e Aplicada, Estr. Dona Castorina 110,
Rio de Janeiro 22460-320, Brazil} \email{senti@impa.br}

\date{\today}

\thanks{Ya. P. was partially supported by the
National Science Foundation grant \#DMS-0503810 and
U.S.-Mexico Collaborative Research grant 0104675. S. S. was supported by a grant from the Swiss National Science Foundation}

\begin{abstract}
For a smooth map $f$ of a compact interval $I$ admitting an
inducing scheme we establish a thermodynamical formalism, i.e.,
describe a class of real-valued potential functions $\varphi$ on
$I$ which admit a unique equilibrium measure $\mu_\varphi$. Our
results apply to unimodal maps corresponding to a positive Lebesgue measure set of parameters in a one-parameter transverse family.
\end{abstract}
\maketitle


\section{Introduction}


In this note we describe some results on thermodynamics for a class of smooth one-dimensional maps. In the classical thermodynamical formalism, given a continuous map $f$ of a compact metric space $I$ and a continuous real valued {\it potential} function $\varphi$ on $I$, one studies the \emph{equilibrium} measures for $\varphi$, i.e., invariant Borel probability measures $\mu$ on $I$ for which the supremum
$$
\sup\Bigl\{h_\mu(f)+\int_I\,\varphi\,d\mu\Bigr\}
$$
is attained, where $h_\mu(f)$ denotes the metric entropy of the map $f$. In this paper we are interested in the existence, uniqueness and ergodic properties of equilibrium measures for a smooth one-dimensional map that admits an ''inducing scheme'' as described below. It represents $f$, restricted to a subset $X\subset I$, as a \emph{tower} $(W, \tau, F)$ where $F$ is the \emph{induced map} acting on the \emph{inducing domain} $W\subset I$ and the \emph{inducing time} $\tau$ is a positive integer-valued function on $W$. The latter is usually not the first return time to $W$. An important feature of the inducing scheme is that $W$ admits a countable generating Bernoulli partition. Thus, $F$ is equivalent to the full shift on a countable set of states.

We apply results of Mauldin and Urba\'nski \cite{MauUrb1} and of
Sarig \cite{Sarig1}, \cite{Sarig2} (see also Aaronson, Denker and
Urbanski \cite{AaDU1}, Yuri \cite{Yuri1} and Buzzi and Sarig
\cite{BuzziSarig1}) to establish the existence and uniqueness of
equilibrium measures for the induced map $F$. We then lift them from the inducing domain to the tower. The latter procedure is quite subtle and makes use of the integrability of the inducing time with respect to the equilibrium measures. It also requires a relation between $F$-invariant measures on $W$ and $f$-invariant measures on $X$ and their respective entropies -- the so called Abramov formula (see Zweim\"uller \cite{Zweimuller7}; for related results see Keller \cite{Keller3}, Bruin \cite{Bru2}).

We give sufficient conditions on the potential function $\varphi$
under which the above procedure can be applied guaranteeing the
existence and uniqueness of equilibrium measures. Let us stress
that we do this without resorting to the study of analyticity of the pressure function (see Remark 3 for more details). Results by
Ruelle \cite{Ruelle2} (see also Aaronson \cite{Aaronson1}) describe some ergodic properties of equilibrium measures for the induced system (including exponential decay of correlations and the Central Limit Theorem) which, by Young \cite{LSY1}, can be transferred to the original system.

We stress that the thermodynamical formalism presented here
depends on the choice of the inducing scheme, since the latter
determines a class of $f$-invariant measures and a class of
potential functions to which our theory applies. Note that one can
construct different inducing schemes for a given map. Naturally one would like the class of measures and potentials, allowed by the scheme, to be as large as possible and, ideally, it should include all $f$-invariant measures and significant potential functions such as constants and $\varphi_t=-t\log|df|$.

In particular, our methods apply to transverse one-parameter families of unimodal maps with parameters from a set of positive Lebesgue measure, where the measures we consider have integrable inducing time and the class of potential functions includes $\varphi_t=-t\log|df|$ with $t$ in an interval containing $[0,1]$ (see \cite{Pesin-Senti2}); this extends results by Bruin and Keller \cite{Bru-Kel} for the parameters and measures under consideration). In particular, this gives a new and unified approach for establishing existence and uniqueness of measures of maximal entropy (first, constructed by Hofbauer \cite{Hof1, Hof2}) as well as absolutely continuous invariant measures. Our results can also be used to observe phase transitions for some potentials associated with uniformly expanding maps of intervals (see \cite{PesinZhang}).

{\bf Acknowledgments.} We would like to thank H. Bruin, J. Buzzi, D.
Dolgopyat, F. Ledrappier, M. Misiurewicz, O. Sarig, M. Viana and M. Yuri for valuable discussions and comments. We also thank the ETH, Z\"urich where part of this work was conducted. Ya. P. wishes to thank the Research Institute for Mathematical Science (RIMS), Kyoto and Erwin Schr\"odinger International Institute for Mathematics (ESI), Vienna -- where a part of this work was carried out -- for hospitality.


\section{Inducing Schemes}\label{inducing}


Let $f:I\to I$ be a smooth map of a compact interval $I$ and $S$ a countable collection of open disjoint intervals. Let also
$\tau:S\to \mathbb{N}$ be a positive integer-valued function.
Define the \emph{inducing domain} by 
$\mathcal{W}:=\bigcup_{J\in S}\,J,$ the \emph{inducing time} 
$\tau\colon \mathcal{W}\to\mathbb{N}$ by
$$
\tau(x):=
\begin{cases}
\tau(J), & x\in J\in S,\\
0 & x\not\in\mathcal{W},
\end{cases}
$$
and  the \emph{induced map} $F: \mathcal{W}\to I$ by
$F(x)=f^{\tau(x)}(x)$.
Denote
$$
W:=\bigcap_{n\geq 0}F^{-n}(\mathcal{W}), \quad
X:=\bigcup_{J\in S}\,\bigcup_{k=0}^{\tau(J)-1}\, f^k(W\cap J).
$$
The set $W$ is totally invariant under $F$ and the set $X$ is forward invariant under $f$. We call $X$ the \emph{tower} with the \emph{base} $W$. In general, $W$ may be an empty set. However, in
many interesting cases including some unimodal maps
(see Section 5) one can construct inducing schemes such that for
an interval $I'$ with $W\subset\mathcal{W}\subset I'\subset I$ the
Lebesgue measure of the set $I'\setminus W$ is zero.

Let $\bar J$ denote the closure of the set $J$. We say that $f$
admits an \emph{inducing scheme} $\{S,\tau\}$ if:
\begin{enumerate}
\item [(H1)] $f^{\tau(J)}|\bar J$ is a homeomorphism onto its image and $f^{\tau(J)}(\bar J)\supseteq\mathcal{W}$;
\item [(H2)] the partition $\mathcal{R}$ of $\mathcal{W}$ induced by the sets $J\in S$ is generating, i.e., for any countable collection of intervals $J_1, J_2, \dots \in S$ the intersection 
${\bar J}_1\cap f^{-\tau(J_1)}({\bar J}_2)\cap f^{-\tau(J_1)-\tau(J_2)}({\bar J}_3)\cap \dots$ consists of a single point.
\end{enumerate}

Condition (H2) often comes as a result of the fact that the induced map is expanding: there exists $\lambda>1$ with $|dF(x)|>\lambda$ for every $x\in\mathcal{W}$.

The partition $\mathcal{R}$ induces a partition of $W$ which we denote by the same letter. Conditions (H1), (H2) allow one to obtain a
symbolic representation of the induced map as the (full) shift $\sigma$ on the space $S^{\mathbb{N}}$ where $S$ is a countable
set of states. More precisely, let
$\mathcal{\tilde{W}}=\bigcup_{J\in S}\bar J$ and define the
\emph{coding map} $h\colon S^{\mathbb{N}}\to\tilde{\mathcal{W}}$ by
$h((a_k)_{k\in \mathbb{N}})=x$ where $x$ is such that $x\in\bar{J}_{a_0}$
and
\[
f^{\tau(J_{a_k})}\circ\cdots\circ f^{\tau(J_{a_0})}(x)\in
\bar{J}_{a_{k+1}} \quad\mbox{ for }\quad k\ge 0.
\]
\begin{proposition}\label{conjugacy}
The map $h$ is well-defined and is onto. It is one-to-one on
$h^{-1}(W)$ and is a topological conjugacy between
$\sigma|h^{-1}(W)$ and $F|W$. Moreover, the set
$S^{\mathbb{N}}\setminus h^{-1}(W)$ is countable.
\end{proposition}
Let $\mathcal{M}(F)$ denote the set of $F$-invariant ergodic Borel probability measures on $W$ and $\mathcal{M}(f)$ the set of $f$-invariant ergodic Borel probability measures on $X$. For any $\nu\in\mathcal{M}(F)$ let
$$\displaystyle{
Q_\nu:=\sum_{J\in S}\tau(J)\nu(J). }
$$ 
If $Q_\nu<\infty$ we define the {\it lifted measure} $\pi(\nu)$ on $I$ in the following way (see for example \cite{Zweimuller7}): for any measurable set $E\subseteq I$
$$
\pi(\nu)(E):=\frac1Q_\nu\sum_{J\in S}\,\sum_{k=0}^{\tau(J)-1}
\nu(f^{-k}(E)\cap J).
$$
Note that $\pi(\nu)\in\mathcal{M}(f)$. If $Q_\nu=\infty$, the measure given by 
$$\sum_{J\in S}\,\sum_{k=0}^{\tau(J)-1}\nu(f^{-k}(E)\cap J),
$$ 
is $\sigma$-finite but not finite.

We call a measure $\mu\in\mathcal{M}(f)$ \emph{liftable} if there is a measure $\nu\in\mathcal{M}(F)$ such that $\pi(\nu)=\mu$. By definition, $Q_\nu<\infty$ and as can be easily seen, $\mu|W\ll~\nu$. It follows that $\nu$ is uniquely defined. We call $\nu$ the \emph{induced measure} for $\mu$ and we write $\nu=i(\mu)$. Finally, we denote the set of liftable measures by $\mathcal{M}^*(f)$.

We also have the following result.

\begin{proposition}[Zweim\"uller \cite{Zweimuller7}]
Any measure $\mu\in\mathcal{M}(f)$ with $\tau\in~L^1(X, \mu)$ is liftable.
\end{proposition}

For $\varphi:I\to\mathbb{R}$ we define the \emph{induced function}
$\tilde\varphi:\mathcal{W}\to\mathbb{R}$ by
$$
\tilde\varphi(x):=\sum_{k=0}^{\tau(J)-1}\varphi(f^k(x))\quad
\text{ for }x\in J.
$$
Observe that given an interval $J\in S$, the function $\tilde\varphi$ can be extended to a continuous function on $\bar J$ which we still denote by $\tilde\varphi$. This extension is well-defined for all $x\in J$, but is "multi-valued" on the (countable) set of points $x\in\overline{J}\cap\overline{J'}$, $J\neq J'\in S$, and the value of $\tilde{\varphi}$ depends on the extension. In what follows the extension that should be used to determine the values of the function at the endpoints of intervals will be clear from the context.
\begin{proposition}[Abramov's and Kac's Formul\ae, see \cite{Zweimuller7}]\label{Abramov}
If $\mu\in\mathcal{M}^*(f)$ then
$$
h_{i(\mu)}(F)=Q_{i(\mu)}\cdot h_\mu(f)<\infty.
$$
If $\int_X\varphi\,d\mu$ is finite then
$$
-\infty<\int_W\tilde\varphi\,d
i(\mu)=Q_{i(\mu)}\cdot\int_X\varphi\,d\mu <\infty.
$$
\end{proposition}


\section{Equilibrium Measures for the Induced Map}


Consider the full shift $\sigma$ on the set $S^\mathbb{N}$, where
$S$ is a countable set of states, and a potential function
$\Phi:S^\mathbb{N}\to\mathbb{R}$.

The \emph{Gurevich pressure} of $\Phi$ is defined by (see for example, \cite{Sarig2}):
\begin{equation}\label{gur}
P_G(\Phi):=\lim_{n\to\infty}\frac1n\log
\sum_{\sigma^n(\omega)=\omega}\exp(\Phi_n(\omega))\chi_{[a]}(\omega),
\end{equation}
where $a\in S$, $\chi_{[a]}$ is the characteristic function of the cylinder $[a]$ and
$$
\Phi_n(\omega)=\sum_{k=0}^{n-1}\Phi(\sigma^k(\omega)).
$$
The \emph{$n$-variation} $V_n(\Phi)$ is defined by
$$
V_n(\Phi):=\sup_{[b_0,\dots,b_{n-1}]}\,\sup_{\omega,\omega'\in
[b_0,\dots,b_{n-1}]}\{|\Phi(\omega)-\Phi(\omega')|\}.
$$
If $\sum_{n\ge 1}\,V_n(\Phi)<\infty$ then the limit in~\eqref{gur}
exists and does not depend on $a\in S$ (see \cite{Sarig1}).

We call a measure $\nu_\Phi$ a \emph{Gibbs measure} for $\Phi$ if
there exist constants $C_1>0$ and $C_2>0$ such that for any
cylinder $[b_0,\dots , b_{n-1}]$ and any $\omega\in [b_0,\dots
,b_{n-1}]$ we have
\[
C_1\le\frac{\nu_\Phi([b_0,\dots ,b_{n-1}])}
{\exp\left(-nP_G(\Phi)+\Phi_{n}(\omega)\right)}\le C_2.
\]
\begin{proposition}[\cite{Sarig2}, see also \cite{MauUrb1},\cite{Aaronson1},
\cite{AaDU1}, \cite{BuzziSarig1}]\label{sarig} Assume that
\begin{itemize}
\item $\Phi$ is continuous and $\sup_{\omega\in S^\mathbb{N}}\Phi<\infty$,
\item $P_G(\Phi)<\infty$,
\item $\sum_{n\ge 1}V_n(\Phi)<\infty.$
\end{itemize}
Then there exists an ergodic $\sigma$-invariant Gibbs measure $\nu_{\Phi}$
for $\Phi$.
\end{proposition}

Observe that a Gibbs measure $\nu_\Phi$ is positive on every
non-empty open set and ergodic and hence, by
Proposition~\ref{conjugacy}, $\nu_\Phi(h^{-1}(W))=~1$.

Given a \emph{potential} function $\varphi: I\to\mathbb{R}$ and its induced potential $\tilde\varphi$, we denote by
$$
\mathcal{M}_{\tilde\varphi}(F)
:=\{\nu\in\mathcal{M}(F)\,:\, -\int_{W}\tilde\varphi\,d\nu<\infty\}.
$$
We call a measure $\nu_{\tilde\varphi}\in\mathcal{M}_{\tilde\varphi}(F)$ an \emph{equilibrium measure} for $\tilde\varphi$ (with respect to the class of measures $\mathcal{M}_{\tilde\varphi}(F)$) if
\[
\sup_{\nu\in\mathcal{M}_{\tilde\varphi}(F)}
\{h_\nu(F)+ \int_W\tilde\varphi\,d\nu \}=
h_{\nu_{\tilde\varphi}}(F)+\int_W\tilde\varphi\,d\nu_{\tilde\varphi}.
\]
Given a potential function $\varphi:I\to\mathbb{R}$, define the
potential function $\Phi:S^\mathbb{N}\to \mathbb{R}$ as follows:
$$
\Phi(\omega):=\tilde\varphi\circ h(\omega)=
\sum_{i=0}^{\tau(J_0)-1}\varphi(f^i(h(\omega)))\quad\mbox{if
}\,\omega=(J_0, J_1, \ldots)\in S^\mathbb{N}
$$
(recall that $\tilde\varphi$ is extended on each $J\in S$ to its closure $\bar J$). Note that $\Phi$ is continuous in the discrete topology of
$S^\mathbb{N}$.

The following result establishes existence and uniqueness of
equilibrium measures for the induced map and a certain class of
potential functions.

\begin{theorem}\label{sarig-induced} Assume that $\Phi$
satisfies the conditions of Proposition~\ref{sarig}. If, in
addition, the entropy $h_{\nu_\Phi}(\sigma)<\infty$, then the
measure
$$
\nu_{\tilde\varphi}:=h_*\nu_\Phi
$$
belongs to $\mathcal{M}_{\tilde\varphi}(F)$ and it is the unique
equilibrium measure for $\tilde{\varphi}$ (with respect to the class of measures $\mathcal{M}_{\tilde\varphi}(F)$).
\end{theorem}
The proof follows from Propositions~\ref{conjugacy} and \ref{sarig} and the fact that $\nu_\Phi(h^{-1}(W))=1$.


\section{Existence, Uniqueness and Ergodic Properties of Equilibrium Measures}


Denote by
$$
s_\varphi:=\sup_{\mu\in\mathcal{M}^*(f)}\left(h_\mu(f)+
\int_X\varphi\,d\mu\right).
$$
A measure $\mu\in\mathcal{M}^*(f)$ for which this supremum is attained is called an \emph{equilibrium measure} for $\varphi$ (with respect to the class of measures $\mathcal{M}^*(f)$).

Our definitions of equilibrium measures differ from the classical
ones. First, we only consider the supremum over measures for which
$W$ (respectively, $X$) is of full measure. Second, we only allow measures $\mu\in\mathcal{M}^*(f)$, i.e., liftable measures. For a general inducing scheme one may not be able to drop these restrictions: Pesin and Zhang \cite{PesinZhang} gave an example of a one-dimensional map $f$ with an inducing scheme and of a potential $\varphi$ for which the supremum over all $f$-invariant ergodic Borel probability measures supported on the closure $\bar X$ of the tower $X$ is strictly bigger than the one taken over $\mu\in\mathcal{M}(f)$ (the former is attained by a measure with $\mu(X)=0$). Also Zhang (following an example of Zweimuller \cite{Zweimuller7}) constructed an example of an abstract tower for which the supremum over all measures $\mu\in\mathcal{M}(f)$ is attained by a measure with 
$\int\tau(x)\,d\mu=\infty$ and is strictly larger than the one taken over measures with integrable inducing time (unpublished).

Given a potential function $\varphi:X\to\mathbb{R}$, we denote  by $\phi^+$ the induced potential of the function $\varphi-s_\varphi$.
Observe that $\phi^+(x)=\tilde\varphi(x)-s_\varphi \tau(x)$.

The class of potential functions $\varphi$ is defined by a collection of axioms that guarantee that the conditions of Theorem~\ref{sarig} hold and that the lifted measure of an equilibrium measure for the induced system is an equilibrium measure for the original system.

More precisely, we assume the following conditions on $\varphi$:

\begin{enumerate}
\item [(P1)] (\emph{boundedness}):
$\displaystyle{\quad \sup_{x\in W}\phi^+(x)<\infty}$;
\item [(P2)] (\emph{local H\"older continuity}): $\varphi$ is continuous on $\mathcal{W}$ and there exist $A>0$ and $0<r<1$ such that $V_n(\tilde\varphi\circ h)\le A r^n$ for all $n\ge 1$;
\item [(P3)] (\emph{finite Gurevich pressure}):
$$
\sum_{J\in S}\,\sup_{x\in J}\,\exp\,\tilde\varphi(x)<\infty;
$$
\item [(P4)] (\emph{positive recurrence}): there exists $\varepsilon_0>0$ such that for every $0\le\varepsilon\le\varepsilon_0$,
$$
\sum_{J\in S}\,\sup_{x\in J}\,\exp\,(\phi^+(x)+\varepsilon \tau(x))<\infty.
$$
\end{enumerate}

The following statement establishes existence and uniqueness
of equilibrium measures for the map $f$.

\begin{theorem}\label{remark}
Let $f$ be a smooth one-dimensional map of a compact interval admitting an inducing scheme $\{S,\tau\}$. Let also $\varphi$ be a potential function satisfying Conditions $(P1)-(P4)$. Then there exists a unique equilibrium measure $\mu_\varphi\in\mathcal{M}^*(f)$ for $\varphi$ (with respect to the class of measures $\mathcal{M}^*(f)$).
\end{theorem}

We outline the proof of the theorem (for the complete proof see
\cite{Pesin-Senti2}). By Conditions (P2) and (P3) respectively, the
induced potential function $\tilde\varphi$ has summable variations and finite Gurevich pressure. This implies that $s_\varphi$ is finite. By Conditions (P1), (P2) and (P4) (with $\varepsilon=0$), the induced potential $\phi^+$ corresponding to the ``normalized'' potential $\varphi-s_\varphi$ is bounded from above, has summable variations and finite Gurevich pressure. Therefore, by Proposition~\ref{sarig}, there exists a unique Gibbs measure $\nu_{\phi^+}:=h_*\nu_{\Phi^+}$ for $\phi^+$ on $W$, where $\Phi^+:=\phi^+\circ h$. Condition (P4) implies that 
$$
\sum_{J\in S}\,\tau(J)\,\sup_{x\in J}\,\exp\,(\phi^+(x))<\infty.
$$ 
One can show that this yields $Q_{\nu_{\phi^+}}<\infty$ and hence,
$-\int_W\,\phi^+\,d\nu_{\phi^+}<\infty$. It follows that
$\nu_{\phi^+}\in\mathcal{M}_{\phi^+}(F)$ and it is an equilibrium measure for $\phi^+$. 
To prove that this  natural candidate is indeed an equilibrium measure for $\varphi$ one needs to verify the \emph{recurrence property}: $P_G(\phi^+)=0$ (see for example \cite{Sarig2}). This can be done using Condition (P4). The variational principle now implies that 
$$
h_{\nu_{\phi^+}}(F)+\int_W\phi^+ \,d\nu_{\phi^+}=0.
$$
By Abramov's and Kac's formulas (see Proposition~\ref{Abramov}), we have that
$$
Q_{\nu_{\phi^+}}\bigl(h_{\pi(\nu_{\phi^+})}(f)+\int_I(\varphi-s_\varphi)\,
d\pi(\nu_{\phi^+})\bigr)=0
$$
since $Q_{\nu_{\phi^+}}<\infty$, and the desired result follows.

In order to study the ergodic properties of the equilibrium measure, we
introduce the following additional condition:
\begin{enumerate}
\item[(P5)]\emph{(exponential tail)}:
there exist $K>0$ and $0<\theta<1$ such that for all $n>0$,
\[
\nu_{\phi^+}(\{x\in W: \tau(x)\ge n\})\le K\theta^n.
\]
\end{enumerate}

\begin{theorem}
Under the same hypothesis as in Theorem \ref{remark} and Condition $(P5)$, the equilibrium measure $\mu_\varphi$ has exponential decay of correlations and satisfies the Central Limit Theorem for the class of functions whose induced function are H\"older continuous on $W$.
\end{theorem}
By Theorem~\ref{remark}, the equilibrium measure $\nu_{\phi^+}$ exists and by results of Ruelle \cite{Ruelle2} (see also Aaronson \cite{Aaronson1}), it has exponential decay of correlations and satisfies the Central Limit Theorem for the class of H\"older continuous potentials on $W$. By results of Young \cite{LSY1}, under Condition $(P5)$, these ergodic properties can be transferred to the equilibrium  measure $\mu_\varphi$.

{\bf Remark 1.}
We call two functions $\varphi$ and $\psi$ {\it cohomologous} if there
exists a bounded function $h$ and a real number $C$ such that
$\varphi-\psi=h\circ f-h+C$. Any equilibrium measure for $\varphi$ is an equilibrium measure for any $\psi$ cohomologous to $\varphi$ and conversely. In particular, if $\varphi$ satisfies Conditions (P1)--(P4) then there is a unique equilibrium measure for any $\psi$ cohomologous to $\varphi$ regardless of whether this function satisfies our Conditions (P1)--(P4).

{\bf Remark 2.} Unlike the classical thermodynamical formalism our approach to establish existence and uniqueness of equilibrium measures does not directly rely on analyticity of the \emph{pressure function} $t\to P_G(\tilde\varphi_1+t\tilde\varphi_2)$ where $\tilde\varphi_1$ and $\tilde\varphi_2$ are the induced potentials for some given potential functions $\varphi_1$ and $\varphi_2$ respectively. On the other hand, using results of Sarig \cite{Sarig3}, one can show that, if the function $\varphi_1$ satisfies Conditions (P2) and (P4), if the function $\varphi_2$ satisfies Condition (P2) and if the function $\varphi_1-s_{\varphi_1}+t\varphi_2$
satisfies Condition (P3) for $t$ near $0$, then $P_G(\varphi_1-s_{\varphi_1}+t\varphi_2)$ is analytic in $t$.


\section{Equilibrium Measures For Unimodal Maps}


We describe applications of our results to families of unimodal maps.


Let $f$ be a $C^3$ interval map with one non-flat critical point
belonging to the interior of the interval of definition of $f$. Without loss of generality, we may assume that $0$ is the critical point and that $f$ is symmetric with respect to $0$, so $f:[-b, b]\to[-b, b]$ for some $b>0$. Assume further that $f(x)=\pm|\theta(x)|^l+f(0)$ for some local $C^1$ diffeomorphism $\theta$ and $l>1$. The map $f$ is \emph{unimodal} if the derivative $df/dx$ changes signs at $0$ and $f(\pm b)\in\{\pm b\}$. Results for a unimodal map $f$ customarily require that $f$ has negative Schwarzian derivative, but this condition can be dropped if $f$ has no attracting periodic points (see
\cite{Graczyk-Sands1}).

A smooth one-parameter family of unimodal maps $f_a$ is called
\emph{transverse} in a neighborhood of a parameter $a^*$ if
$$
\frac{d}{da}f_a(0)\ne\frac{d}{da}\zeta(a),
$$
where $\zeta(a)$ is the smooth continuation of the point
$x^*:=\zeta(a^*):=f_{a^*}(0)$. For any such family, Yoccoz
\cite{JCY3} (see also \cite{Senti3}) introduced the set
$\mathcal{A}$ of \emph{strongly regular} parameters, which has
positive Lebesgue measure (such parameters also satisfy the
Collet-Eckmann condition).

Fix $a\in\mathcal{A}$ and set $f=f_a$. Let
$I_{[b_0,\dots,b_{n-1}]}$ denote the set of points $x$ with
$F^i(x)\in I_{[b_i]}\in S$ for $0\le i\le n-1$ where $F$ denotes
the induced map for $f$.

\begin{theorem}[\cite{JCY3},\cite{ Senti3},\cite{Pesin-Senti2}]
Consider a transverse one-parameter family of unimodal maps with
no attracting periodic points. There exists a positive Lebesgue measure set $\mathcal{A}$ of parameters such that for every $a\in\mathcal{A}$ the map
$f=~f_a$ admits an inducing scheme $\{S,\tau\}$ satisfying
(H1), (H2). In addition, there exists an interval $I'$ with
$W\subset I' \subset I$ such that Leb$(I'\setminus W)=0$ and the
following conditions hold:
\begin{enumerate}
\item [(H3)] there are constants $c_1>0$ and  $\lambda_1>1$ such that for all $n\ge 0$,
$$
\sum_{J\in S: \,\tau(J)\ge n}|J|\le c_1^{-1}\lambda_1^{-n};
$$
\item [(H4)] (Bounded Distortion): for each $n\ge 0$, each interval
$I_{[b_0,\ldots,b_{n-1}]}$ and each $x, y\in I_{[b_0,\ldots,b_{n-1}]}$,
$$
\left|\frac{dF^n(x)}{dF^n(y)}-1\right|\le c_2 |F^n(x)-F^n(y)|;
$$
\item[(H5)] for every $\gamma>1$ there exists $c=c(\gamma)$ such that
$$
Card\{J\in S\ |\ \tau(J)=n\}\le c\gamma^n.
$$
\end{enumerate}
\end{theorem}
As a corollary of this theorem one has the following estimates: there exist positive constants $c_2,c_3, c_4$ and $\lambda_2\ge\lambda_1$ such that for every $J\in S$ and $x\in J$,
\begin{equation}\label{bd}
c_1\lambda_1^{\tau(J)}\le |J|^{-1}\le c_2\lambda_2^{\tau(J)}, \quad
c_3|J|^{-1}\le |dF(x)|\le c_4|J|^{-1}.
\end{equation}
Given a map $f_a$ from a transverse family of unimodal  maps, we consider the potential function
$$
\varphi_{a,t}(x)=-t\log|df_a(x)|.
$$
If the parameter $a$ is strongly regular, using (H3)--(H5) and \eqref{bd}, one can show that the function $\varphi_{a,t}$ satisfies Conditions (P1)--(P4) and hence, results of the previous section apply.
\begin{theorem}[\cite{Pesin-Senti2}]
Consider a transverse one-parameter family of unimodal maps with
no attracting periodic points. There exists a positive Lebesgue measure set $\mathcal{A}$ of parameters such that for every $a\in\mathcal{A}$ the following
statements hold.
\begin{enumerate}
\item There exist $t_0(a)<0<1<t_1(a)$ such that for every $t_0(a)<t<t_1(a)$ there is a measure $\mu_{a,t}$ on $I$ which is the unique equilibrium measure for $\varphi_{a,t}$ (with respect to the class of measures $\mathcal{M}^*(f)$). Moreover, for any measure $\mu$ not supported on $X$, 
$$
h_\mu(f_a)-t\int\log |df_a(x)| d\mu
< h_{\mu_{a,t}}(f_a)-t\int\log |df_a(x)| d\mu_{a,t}.
$$
\item The measure $\mu_{a,t}$ is ergodic, has exponential decay of
correlations and satisfies the Central Limit Theorem for the class of
functions whose induced functions are H\"older continuous. The measure $\mu_{a,1}$ is the absolutely continuous invariant measure and the measure $\mu_{a,0}$ is the unique measure of maximal entropy.
\end{enumerate}.
\end{theorem}

\bibliographystyle{alpha}
\bibliography{MMJ2005}
\end{document}